# Demonstrative and non-demonstrative reasoning by analogy

Emiliano Ippoliti

Analogy and analogical reasoning have more and more become an important subject of inquiry in logic and philosophy of science, especially in virtue of its fruitfulness and variety: in fact analogy «may occur in many contexts, serve many purposes, and take on many forms»[1]. Moreover, analogy compels us to consider material aspects and dependent-on-domain kinds of reasoning that are at the basis of the lack of well-established and accepted theory of analogy: a standard theory of analogy, in the sense of a theory as classical logic, is therefore an impossible target. However, a detailed discussion of these aspects is not the aim of this paper and I will focus only on a small set of questions related to analogy and analogical reasoning, namely:
1) The problem of analogy and its *duplicity*;
2) The role of analogy in demonstrative reasoning;
3) The role of analogy in non-demonstrative reasoning;
4) The limits of analogy;
5) The convergence, particularly in *multiple* analogical reasoning, of these two apparently distinct aspects and its philosophical and methodological consequences;

## § 1 The problem of analogy and its duplicity: the controversial nature of analogical reasoning

One of the most interesting aspects of analogical reasoning is its intrinsic *duplicity*: in fact analogy belongs to, and takes part in, both demonstrative and non-demonstrative processes. That is, it can be used respectively as a means to prove and justify knowledge (e.g. in automated theorem proving or in confirmation patterns of plausible inference), and as a means to obtain new knowledge, (i.e. in heuristics and in the generation of conjectures and hypotheses). As it is well known, these two kinds of reasoning are traditionally (e.g. by the logical empiricist philosophy of science) treated as distinct and belonging to two completely independent contexts, namely the context of justification on one hand and the context of discovery on the other: *justification* is the phase in which hypotheses are confirmed or rejected, *discovery* is the phase in which the scientific hypotheses are generated. Moreover, analogical reasoning is widely considered not only as one of the main tools in problem-solving activity, but also as an ubiquitous, highly controversial and complex concept: in fact «it can be said that the analogy is, as the tongue of Aesop, at the same time the best and the worse thing»[2]. Its controversial nature is not accidental and relies on two fundamental properties of analogical reasoning:

---
[1] Shelley 2003, 1
[2] Dieudonné 1981, 257

a) *Ampliativity* (i.e. the capability to really extend our knowledge by reaching conclusions which are not included in the premises).

b) *Non-monotonicity* (i.e. the sensitivity to the entry of new information and premises, which are able to modify the previously obtained conclusions). As a consequence, analogy is an intrinsic time-sensitive kind of inference: it strictly depends on the background knowledge existing at a given time.

Moreover, the very definition of analogy is problematic. The existing orthodox literature agrees in considering analogy as a kind of comparison, which, in short, allows to transfer a known property/information from a sufficiently known *source* domain S to an at least partially unknown *target* domain T, by a relation of mapping $\mu$ of objects, relations and properties from S into T. In particular, it is possible to distinguish two main conceptions of analogy, namely the *inductive* and the *structural*.

1) **Analogy as induction (*inductive conception*)**
Analogy is a form of induction (induction on attributes or properties), in virtue of which a single observation is used as a basis for some conclusion. In this sense analogy is a kind of generalization (e.g. Keynes[3]), which is obtained by a conjunction of material resemblances between domains.

2) **Analogy as shared structure (*structural conception*)**
Analogy is a mapping or alignment of «hierarchically structured, causal relationships shared between source and target analogs»[4]. That is, analogy is ideally an isomorphism of two domains (e.g. Hempel's nomic isomorphism between Ohm's law and Poiseulle's law).

The two conceptions agree on the relevance of overall similarity across domains; nevertheless, the structural conception is based on the mapping between relations (and not on attributes as the inductive conception) and on the *systematicity principle*, which claims that an analogy is *good* if it contains mapping between richly structured higher order relations (which are in general the causal ones). Therefore, both conceptions try to specify «a rationale for analogical reasoning»[5], that is, to offer an answer to LPA, the Logical Problem of Analogy.

LPA can be formulated as the problem to «find a criterion which, if satisfied by any particular analogical inference, sufficiently establishes the truth of the projected conclusion for the target»[6]. Indeed, LPA deals with the question of the *goodness* of an analogy. Obviously, the truth or the falsity of the conclusion of an analogical inference cannot be a way of evaluating the inference and do not say anything about its goodness: in fact, «the validity of empirical arguments as the foundation of a probability cannot be affected by the actual truth or falsity of their conclusions»[7]. As a matter of fact, analogical inference is context-dependent and time-sensitive: an analogical projection which is false at time $t_1$ and within a certain context, can become true at time $t_2$ and within another context (e.g. extended background knowledge). Therefore, there is a set of traditional criteria that can be used to evaluate the goodness of analogies and to give an answer LPA.

---

[3] Keynes 1921
[4] Shelley 2003, 7
[5] Leblanc 1969, 29
[6] Davies 1988, 229

For example[8], it is possible to point out the following seven common criteria: (1) more similarities between source and target imply a stronger analogical argument; (2) more differences imply a weaker argument; (3) a greater ignorance about the items compared imply a weaker argument; (4) a weaker conclusion implies more plausible argument; (5) analogical arguments based on causal relations tend to be more plausible then those which are not based on causal relations; (6) structural analogies are stronger than those based on superficial similarities; (7) the relevance to the conclusion of similarities and differences should be considered. Therefore the theory of similarity, the theory of relevance and determinant structures, and the theory of typicality can be used to deal with LPA: their adoption is strictly connected to what conception of analogical reasoning is chosen.

Obviously, all these theories fail in solving LPA, which, in line of principle, is not solvable because it is subject to the *paradox of inference*[9]: you cannot have a genuine analogy which is both sound and ampliative. Therefore, if you want to have a sound analogical inference, you have to renounce to its very nature, the ampliativity, i.e. you have to transform it into a valid argument. But such a transformation is possible only in trivial cases of analogical transfer and under very restrictive, accurate and uninteresting conditions which, definitively, make the transfer not really analogical.

§ 1.2  **Theory of similarity**

The theory of similarity is the main model developed by the inductive approach and its theorization can be traced back to the very origin of theory of analogy, i.e. to Mill's competition of similarities and dissimilarities (in Keynes's terminology, positive and negative analogies, that is shared features and differences respectively). According to the conception of analogy based on similarity «the general form of analogical inference is: if *P(a)* and *b* è similar to *a*, then *P(b)*. This form presupposes that the meaning of 'similar' has been specified. Now, there are several notion of similarity»[10]. This conception tries to solve LPA by elaborating a "*metrics*" of analogy, that is, by the refinement of the concept of *degree* of similarity. Similarity-based analogical inference is usually defined in terms of aspects[11]: it allows inferring new similarities between the source S and the target T *w.r.t. an aspect* Q from known similarities between S and T *w.r.t. an aspect* P.

---

[7] Keynes 1921, 270
[8] Bartha 2002
[9] Cellucci 2002
[10] Cellucci 1998, 367
[11] Melis - Veloso 1998

$$\frac{Similar\ (S,T,P)}{Similar\ (S,T,Q),\ Q(S)}$$
$$Probably\ Q(T)$$

An aspect is therefore defined as a part of the object's description (e.g. a set of feature-value pairs). It is worth noting that such a definition is quite problematic: although the representation of both source and target in a formal language is common in analogical reasoning, mixed representations of knowledge (e.g. visual representation) are becoming more and more important in the literature on the subject.

The rationale for the use of similarity as a basis for the evaluation of analogical reasoning is a pure and rather intuitive probabilistic argument: if two objects resemblance each other under one or more aspects, then *probably* they will match under other aspects. In fact, let suppose that the two objects $O_1$ and $O_2$ match for the value *v* under *i* of the *m* aspects $A_1$, $A_2$, ... , $A_m$ and let *pr(s,l)* express the probability that $O_1$ and $O_2$ will match under successive aspects $l = i - m$. Then we have that:

$$(G)\ pr(s,l) = \frac{i-m}{m}$$

That is, under conditions of ignorance (i.e. posing all the aspect as having the same probability), (G) expresses numerically the *degree of similarity* between $O_1$ and $O_2$. This equation simply claims that the higher the number of similarities (the matches on values of attributes), the higher the probability that $O_1$ and $O_2$ will match on successive attributes (posing *l* as finite). This simple probabilistic ground allows to consider the analogical inference based on similarity as a kind of induction (a generalization), where:

«the whole process of strengthening the argument in favour of generalization $g(\varphi, f)$ by the accumulation of further experience appears to me to consist in making the argument approximate as nearly as possible to the conditions of a perfect analogy, by steadily reducing the comprehensiveness of those resemblances $\varphi_l$ between the instances which our generalization disregards»[12].

Different measurers of similarity (expressible by indexes) can be used. A set of them individuates the *geometrical* model of similarity[13], which is defined as a distance *D(i,j)* on an opportune *metrics*, according to which the closer two objects *i, j* are in the metrics, the more similar they will be (their comparison is numerical). Usually, these indexes require two steps: at first, *local* similarity – *sim* – is set between the characteristics of a set of characteristics and then combined in the construction of a *global* similarity – *SIM* –. For example[14], it is possible to set local indexes of similarity between features from a given set of features as follows:

i) $$sim(a,b) = \frac{card(a \cup b) - card(a \cap b)}{card(a \cup b)}$$

where *card* is the size of the set, and *a,b* are sets of feature values.

---

[12] Keynes 1921, 253
[13] Torgerson 1977

ii) $$sim(a,b) = \begin{cases} 0 : a \cap b = \emptyset \\ 1 : a \cap b \neq \emptyset \end{cases}$$

"primarily for sets of symbolic feature value"[15].

iii) $$sim(a,b) = \frac{|a-b|}{ul}$$

for «numeric feature values, where |a-b| is the absolute value of the difference between *a* and *b*, and *ul* is the absolute value between the upper and the lower bounds of the interval of possible values of the features»[16]. Once defined over every feature, these local measures are combined in global measures (indexes) of similarity between two objects A and B. Standard indexes of global similarity *SIM* are set, e.g., as follows[17]:

i) *City-Block* similarity $$SIM(A,B) = \frac{1}{p} \sum_{i=1}^{p} sim_i(a_i, b_i)$$

where *p* is the number of the characteristics, $sim_i$ is the local similarity of *i*-th feature $a_i$ and $b_i$.

ii) *Euclidean* similarity $$SIM(A,B) = \sqrt[2]{\frac{1}{p} \sum_{i=1}^{p} (sim_i(a_i, b_i))^2}$$

iii) *Simple Matching Coefficient* $$\frac{\alpha + \delta}{\alpha + \beta + \gamma + \delta}$$

where: α = number of agreeing, positive features; β = number of inconsistent but positive features; γ = number of inconsistent negative features; δ = number of agreeing negative features.

---

[14] Melis - Veloso 1998; Goldstone 1999
[15] Melis - Veloso 1998, 15
[16] *ibid.*
[17] *idib.*

Now, the different indexes better fit different contexts, but they all satisfy the following three metrical properties (and in this sense they are "*geometrical*"):

- *Symmetry*: $SIM(A,B) = SIM(B,A)$
- *Triangle inequality*: $SIM(A,B) + SIM(B,C) \geq SIM(A,C)$
- *Minimality*: $SIM(A,B) \geq SIM(A,A)=0$

These properties are obviously too restrictive and all empirically easily violable, as pointed out by Tversky[18]. For this reason, Tversky develops an alternative and weaker approach, known as *Contrast Model*. It is based on measuring matching features of compared entities combined by the formula:

$$SIM(A,B) = \alpha f(A \cap B) - \beta f(A-B) - \gamma f(B-A)$$

where $(A \cap B)$ expresses the number of common features between entities A and B, (A-B) represents the features that A has but does not B, and (B-A), represents the features that B has but A does not. $\alpha$-$\beta$-$\gamma$ are weights on common and distinctive components and *f* is a function (commonly assumed to be simply additive).

Another classical measure of similarity is the *transformational distance*[19], according to which the similarity of two entities is defined as inversely proportional to the number of legal operations required to transform one entity so as to be identical to the other. For example, let us suppose to have the three following strings of two symbols (+,-):

a) +++---
b) +++--+
c) +--+++

Moreover, let us suppose that the allowed operations are substitution and inversion: thus, while the string (a) requires only just one operation to be transformed so as to be identical to (b) – the simply substitution of "-" with "+" in the 6<sup>th</sup> position from the left, the string (c) requires two operations – the substitution of "–" with "+" in the first position from the left and then the inversion of the string. Therefore, the string (b) is more similar to (a) than the string (c) w.r.t. transformational distance.

The aim of this paper is not to give the list of all indexes of similarity and their several refinements. All these indexes try to account for different aspects of similarity between objects (entities), and define different theories of similarity. Unfortunately, the inductive theories and their measures not only do not solve LPA (i.e. an arbitrary high degree of similarity does not justify at all the adoption of an hypothesis or conclusion based on an analogy: the analogical jump can be merely accidental), but they are also affected by a deep and unavoidable weakness, an inner limit: they do not consider the *relevance* of similarities and dissimilarities between the compared objects. In fact, even though the degree of similarity is arbitrarily high, the analogical inference can be completely unjustified because it misses also only one *relevant* feature. Therefore, similarity between objects is not able to say anything about the projected conclusions: the analogical transfer is always an audacious and hazardous passage, in which the reach for ampliativity implies, as required by the paradox of inference, the renounce to every kind of justification, both logical and probabilistic. In particular in the second case, «neither the classical, nor the modern theory of probability have been able to give a satisfying description and justification of inference by analogy»[20], and, I would add, this is not a lack of the probability theory, but a characteristic of the intrinsic nature of analogy, that cannot be really "well-founded" by similarity or, as I am going to show, by other, alternative concepts.

---

[18] Tversky 1977
[19] Hahn - Chater 1997

### § 1.3 Theory of relevance

Theory of relevance is the main device used by the structural approach. It explicitly tries to overcome the inner limit of similarity by taking into account the concept of relevance, that is the relevance of the property for the analogical jump.

In fact theory of relevance seems to offer a reliable criterion for individuating conditions that justify an analogical inference: «these conditions belong to the background knowledge used in analogical reasoning and, hence, justified analogical inference is considered knowledge-based opposed to blind similarity-based analogical inference»[21].

In particular, this criterion is represented by the concept of connection (or the stronger determination): in fact «connections justify analogical inference at some extent»[22]. Connection is a kind of knowledge which formalizes a sort of dependence of an aspect Q from another aspect P, denoted, e.g.[23], by [P,Q]. It is worth noting that connections are weaker then implications: in fact P *determines* Q $\leftrightarrow$ $\exists R \forall \overline{x}\overline{y}\overline{z}(P(\overline{x},\overline{y}) \rightarrow (Q(\overline{x},\overline{z}) \leftrightarrow R(\overline{y},\overline{z})))$, for a relation R which match values of R with values of Q, while P *implies* Q $\leftrightarrow$ $\forall \overline{x}\overline{y}\overline{z}(P(\overline{x},\overline{y}) \rightarrow Q(\overline{x},\overline{z}))$.

Therefore, the structural approach focuses the attention only on relevant or determinative properties (*determinant structures*[24]) during the construction of the analogical jump. The rule for inferring by analogy under determination can be expressed in two main ways, i.e. the complete or strong form (1) and the incomplete or weak form (2), on the basis of the kind of connections:

(1)     **DET$_1$**

---

[20] Carnap 1950, 569-570  
[21] Melis - Veloso 1998, 17  
[22] *ibid.*  
[23] *ibid.*  
[24] Weitzenfeld 1984

$$\frac{P(S) = P(T), [P,Q]}{Q(S) = Q(T), Q(S)}$$
$$Q(T) \ (t)$$

(2) **DET$_2$**

$$\frac{P(S) = P(T), [P,Q]}{\textit{probably}( Q(S) = Q(T)), Q(S)}$$
$$\textit{probably} \ (Q(T) \ (t))$$

These rules express a process according to which once you have a connection [P,Q], you can go on to search a source *S* such that P(S) = P(T), where T is the target, and for which is probable that G(S) = G(T); then, from the knowledge of G(S) (s), you can project G(S) as an approximate value of G(T) (t). Therefore determination theory offers an alternative criterion for evaluating an analogy and solving LPA, in which apparently «the overall degree of similarity does not play any role in the process»[25].

As I showed, in its stronger version the structural approach simply claims that a determination is able to justify the analogical inference and to solve LPA, that is to make an analogical inference well-founded, offering a criterion for «valid reasoning by analogy»[26].

Now, the real question is: does the structural approach really solve LPA?
In the first case, we define as *total* the connections in which P includes all the pieces of information relevant to Q: therefore, the connection [P,Q] completely justifies the analogy. But, it can be noted that it is not an analogical inference at all: it is simply a deductive inference, maybe able to «provide abbreviations of deductive reasoning procedures and thus satisfy one of the advantageous features of analogy: increasing efficiency»[27]. Therefore, it is sound, but not ampliative.

In the second case, we simply define the connections that are not total as *incomplete*, i.e. where P does not include all the pieces of information relevant for Q. Again they do not solve LPA: just like in the case of analogical inferences based on similarity, incomplete-connections-based analogical inferences are not justified at all (they could be accidental). As a matter of fact, «analogical inferences by incomplete connections [P,Q] commonly require further modifications of the result»[28]: they link premises and conclusions in a not stronger way than by similarity. Therefore, they are blind.

Ultimately, it is worth noting that similarity, even though it is explicitly contested by the theory of relevance, can somehow play an important role in it. In fact, a probabilistic and inductive argument can be used to choose the better relevant source for analogical transfers: under the condition of ignorance, and then assuming equally likelihood, the higher the number of similarities, the greater the chance that source will match on relevant attributes is. More formally, let us suppose that S and T match on the values of *s* of the aspects $A_1, A_2, ... , A_m$; moreover let us suppose that the *j* of the *m* aspects are relevant for Q – the property projected form the source to the target. Let *pr*(s,*j*) be the probability that S will match T under the relevant attributes *j*. Then we have that

$$pr(s,j) = \frac{\binom{s}{j}}{\binom{m}{j}}$$

i.e. the higher the number of similarities, the greater the possibility of match on relevant features is.

---

[25] Russell 1988, 252
[26] *idib.*
[27] Melis - Veloso 1998, 19
[28] *idid.*

## § 1.4 Theory of typicality

Analogical inference based on typical instances (instances that are «representative of a certain situation or concept»[29], i.e. having the minimal number of features common, or transferable, to each instance of the same situation or concept) is, especially under certain conditions, another natural candidate to solve LPA. In fact «psychological investigations provided evidence that typical instances are preferred sources for an analogical reasoning where source and target are instances of the same category»[30]. Theory of typicality can be seen both as an inductive and as a structural model of analogy. In other words, typicality can be conceived both as a special case of determination, and therefore able to justify to some extent the analogical inference, and as a kind of similarity. On one hand, it is a kind of generalization (induction), i.e. generalization from a single case (an exemplar case). On the other hand, it includes relevance, even though a slightly different kind of relevance: in fact «the relevance information used here differs from the relevance knowledge assumed for the connections-based analogical inferences»[31]. In fact, a feature Q is relevant for the typical instance S when it can be legitimately hold as common (or at least transferable) to every other instances of the considered concept/situation. Thus, for example, «if Berlin is a typical city and its transportation system include an underground, buses, and taxis and if Berlin is the typical example for Rome in the categories of cities (i.e. Rome is similar to Berlin), then Rome's transportation system includes an underground, buses and taxis»[32].
The conceptual structure underpinning the notion of typicality is then expressible as an acyclic oriented graph: each instance within this conceptual structure is ordered on the basis of its features, according to a relation that is reflexive, transitive and anti-symmetric. More formally, it is possible to conceive typicality «as a partial order $\sqsubseteq$ over a set of instances of a concept C»[33] that provides an evaluation, or an index, of typicality.

---

[29] Kerber - Melis 1997, 157
[30] Melis - Veloso 1998, 20
[31] *idib.*
[32] ivi, 21

Now, in particular, theory of typicality adopts a definition of a concept as a «couple ($\varepsilon_k, \sqsubseteq_k$) where $\varepsilon_k$ is a set of examples and $\sqsubseteq_k$ a partial order on $\varepsilon_k$ (the evaluation or index of typicality)»[34] and an exception is «an example that is not comparable with any other example of $\varepsilon_k$ (i.e., it doesn't exist a e'$\in\varepsilon_k$ with e $\sqsubseteq_k$ e' or e' $\sqsubseteq_k$ e»[35]. A typical example *e* of a category k is now definable as a «maximal (i.e. it doesn't exist any e'$\sqsubseteq_k$ e), which isn't an exception»[36].

As a consequence, the rule, say TYP, to "formalize analogical inference based on typical instances uses a context-independent typicality rating $\sqsubseteq_k$ [...] as a partial order relation in concept structure C and the relevance knowledge "Q is relevant for the typical instance S"»[37]:

**TYP**

$$\frac{tipex(S),\ T\sqsubseteq_k S,\ rilevant(Q,S)}{Q(S) = Q(T) \qquad (QS)}$$
$$Q(T)\ (t)$$

As the DET rules, also TYP rule expresses an exact process, according to which once you have a typical instance S (the source), such that T (the target) is T$\sqsubseteq_k$ S and you know that G is relevant for S, then you can project the value of source G(S) on the target, obtaining G(T).

Unfortunately, also typical-instance-based analogical inference is not really able to offer a relevant contribute to solve LPA.

First of all, typicality is a controversial concept: if we accept its description as partial order relation over a set of instances of a concept, then it is possible to easily show that there are several typical instances for the same category[38]. This complicates a lot the typical-instance-based analogical inference.

Secondly, typicality is affected by the same limits of relevance: the background knowledge added in order to legitimate the analogical jump, makes it, definitively, redundant. The knowledge obtained by its application is *de facto* encapsulated in its category: therefore you have a transfer which is sound, but which is not analogical at all (it is not ampliative).

---

[33] Kerber - Melis - Siekmann 1992, 8
[34] ivi,10
[35] *ibib.*
[36] *ibib.*
[37] Melis - Veloso 1998, 20
[38] Kerber - Melis 1997

Moreover, typicality can be conceived as expressing a distance, i.e. a kind of similarity: in fact the typical instance is the one having the minimal set of features, common to each instance of the concept taken into consideration (i.e. it is a kind of average distance). But typicality isn't justified at all: it shows only a kind a similarity between instances, which can be merely accidental and not founded at all.

So again, typicality can be useful, al least, to «provide abbreviations of deductive reasoning procedures and thus satisfies one of the advantageous features of analogy: increasing efficiency»[39]. Therefore, if it is conceived as sound, it cannot be ampliative, and vice versa. Thus LPA can't be solved by means as similarity, determination or typicality.

## § 2 Demonstrative reasoning by analogy

Analogical reasoning plays an important role in demonstrative reasoning (i.e. as means of justification), in particular in theorem proving and in processes of corroboration of conjectures and hypotheses:

A) In theorem proving the analogical relation among proofs is a powerful means to produce new proofs from existing ones (i.e. generation of proofs by analogical transfer of *proof-outline*[40] or *proof-plan*[41]).

B) In the justification of hypotheses, analogy plays a decisive role in finding analogous sources in order to increase their credibility.

A) Various empirical studies[42] shows that analogical reasoning can be a tool to generate new proofs from existing proofs: «methods, rather than just single calculus steps, are transferred analogically in mathematical theorem proving by analogy»[43], where a method is «a (partial) specification of a tactic, represented in a meta-language, where a tactic executes a number of logical inferences»[44]. In particular, in automated theorem proving, analogy is useful for long and complex proofs: it allows the transposition of a long, complex proof-plan or a proof-outline of some mathematical statements related to one particular domain of mathematics to other statements, either in the same domain or related to another domain. Proof-plan and proof-outline represent the structure of a proof and can be here considered as equivalent (it is worth noting that the search for a *clear* proof-outline or proof-plan is an essential part of mathematical practice).

---

[39] Melis - Veloso 1998, 19
[40] Van Bendegem, 2000
[41] Melis, 1995
[42] e.g. Melis, 1993; Melis, 1994; Bendegem, van 2000
[43] Melis - Veloso 1998, 54
[44] *ibid.*

The transfer of a proof-structure can take many forms: in fact «many proofs by analogy result from transferring parts of the source proof to parts of the target proof, and some proofs by analogy transfer only the proof idea, while others transfer a detailed proof»[45]. Now «in many of these analogies, symbol mapping is insufficient. A change of a problem representation by, e.g., unfolding a definition, can be necessary to reveal the commonality of theorems or assumptions, upon which the analogical transfer is based»[46]. The typical use of analogy in theorem proving is the search of an extension of dimension, i.e. a kind of generalization, for example a proof-outline of a statement about objects (or an object) in two dimensions into a proof-outline of objects (or one object) in three (or more) dimensions. A nice example[47] of proof-plan transfer by analogy is the transfer of Heine-Borel theorem (HBT) from dimension one to dimension two, i.e from $R^1$ to $R^2$.

B)  Analogy can play an essential role in the corroboration of hypothesis, by offering a method for searching for analogous cases in order to reinforce the hypothesis in question. The analysis of the process of corroboration of Euler's solution of the problem of Megoli offered by Polya is a straightforward example of justification by analogy (see §3). Polya[48], in fact, gives to analogy a central role in his patterns of plausible inferences, i.e. pattern of corroboration of hypothesis[49].

## § 3    Non-demonstrative reasoning by analogy

Even if analogy plays an important role in demonstrative reasoning, it is in non-demonstrative reasoning (i.e. as a means of formulating conjectures and hypotheses) that it reveals all its power and centrality in reasoning, also in mathematics, which is rich of examples of analogy as a heuristics means. For example it exists a well-established literature that recognizes the decisive role of analogy in non-demonstrative reasoning such as in:

- knot theory (analogy between knots and numbers) and its development[50]
- the birth of algebraic topology[51]
- the discovery of quaternions by Hamilton[52]
- Euler solution of the *problem of Mengoli* and the Leibniz's series[53]

---

[45] Melis - Veloso 1998, 55
[46] Melis - Veloso 1998, 54
[47] Melis 1995
[48] Polya 1954, I
[49] Ippoliti 2005
[50] Brown 1999, Deninger 1998
[51] Sikora 2001
[52] Van De Waerden 1976
[53] Polya 1954

Euler's solution of Mengoli problem (say MP) is a classic topic in studies dedicated to the analysis of the role of analogy in demonstrative and non-demonstrative reasoning, both in mathematics and natural sciences[54].

I will discuss only this example because I don't totally agree with the current analysis of it. In fact, many authors rely on a standard account of analogy and do not recognize that a more accurate analysis of it can suggest a third (beyond inductive and structural) and more realistic conception of analogical reasoning, i.e. the *heuristics conception*, which accounts for analogy considering the paradox of inference as constitutive of it. According to the heuristic conception, the main question is not to solve LPA, but to offer methods to solve problems.

The problem of Mengoli concerns the determination of the value of the series:

(MP) $1 + \frac{1}{4} + \frac{1}{9} + \frac{1}{16} + \ldots + \frac{1}{n^2} = ?$

The process followed by Euler to find the solution can be summarized in the following way. It consists in six passages, which can differ essentially from its current accounts.

1) In order to solve the problem (MP), Euler tries to go back to a known analogous result, that satisfies some of the properties which constitute the conditions of solvability of the problem given by its preliminary analysis: the question whether the problem of Mengoli is solvable or not is then reduced to the one of finding the value of a known result, or a combination of known results, that share such properties. The list of these properties, given by the analysis of the conditions of solvability, is the following:

   (a)   let it be expressible as an infinite series
   (b)   let it be expressible as the following fraction:
   $\frac{1}{x_1^2} + \frac{1}{x_2^2} + \frac{1}{x_3^2} + \ldots + \frac{1}{x_n^2}$
   (c)   let it be $x_1 = 1y$, $x_2 = 2y$, ..., $x_n = ny$

   Thus, the target of the problem becomes the search for a known-valued series which satisfies (a) and (b) and (c).

2) At this point, the *interaction* with the corpus of existing knowledge at Euler's time (i.e. *Algebra*) allows taking into account an algebraic equation that satisfies the property (b) - a *positive analogy* - but not the property (a) (it isn't an infinite series) - a *negative analogy*. Such an equation is (d):

---
[54] Polya 1954, Van Bendegem 2000, Bartha 2002, Corfield 2003

$$b_1 = b_0 \left( \frac{1}{\beta_1^2} + \frac{1}{\beta_2^2} + \frac{1}{\beta_3^2} + \ldots + \frac{1}{\beta_n^2} \right)$$

which arises from (e)

$$b_0 \left(1 - \frac{x^2}{\beta_1^2}\right)\left(1 - \frac{x^2}{\beta_2^2}\right) \ldots \left(1 - \frac{x^2}{\beta_n^2}\right)$$

which, for $b_0 \neq 0$, expresses the relation between the coefficients and the roots of a generic algebraic equation (f)

$$b_0 - b_1 x^2 + b_2 x^4 - \ldots + (-1)^n b_n x^{2n+1}$$

which has (g) $2n$ many roots $\beta_1, -\beta_1, \beta_2, -\beta_2, \ldots, \beta_n, -\beta_n$.

3) At this point Euler makes his first "analogical jump": in fact he formulates by analogy the hypothesis (I) that a property, i.e. (d), which holds for finite cases, can be valid for infinite cases too. So the initial problem (MP) is reduced to the search for a (d)-type of function which satisfies the conditions (f) - (g) - (c).

4) A new interaction with the corpus of existing knowledge available at Euler's time (i.e. *Trigonometry*) suggests to the Swiss mathematician that there exists an infinite series (c) which can be expressed as (f): it's the development in power series of (h) $sin(x) = 0$, which is equal to

$$\frac{x}{1!} - \frac{x^3}{3!} + \frac{x^5}{5!} - \frac{x^7}{7!} + \ldots = 0$$

which has $2n+1$ many roots: $0, \pi, -\pi, 2\pi, -2\pi, \ldots, n\pi, -n\pi$.

5) Now it is analytically possible to transform the development by means of powers series of *sin(x)* into an infinite algebraic equation which is analogous to (b) by simply dividing the two members of the above equation by $x^1$, that is the linear factor which corresponds to the root 0:

$$\frac{\sin x}{x} = 1 - \frac{x^2}{3!} + \frac{x^4}{5!} - \frac{x^6}{7!} + \ldots$$

That is, the (f)-type equation with $2n$ many roots $\pi, -\pi, 2\pi, -2\pi, \ldots, n\pi, -n\pi$ (c).

6) Now, by analogy, (i.e. on the basis of the matched similarities), Euler makes a second "jump": he conjectures in fact the hypothesis (II) that properties - i.e. (e) and (d) - holding for algebraic equations, will hold for non-algebraic equation (trigonometric) too. Therefore, he represents $\frac{\sin x}{x} = 0$, according to (e), as the infinite product

$$\left(1 - \frac{x^2}{\pi^2}\right)\left(1 - \frac{x^2}{4\pi^2}\right)\left(1 - \frac{x^2}{9\pi^2}\right)\ldots$$

that is, according to (d),

$$\frac{1}{3!} = \frac{1}{\pi^2} + \frac{1}{4\pi^2} + \frac{1}{9\pi^2} + \frac{1}{16\pi^2}\ldots$$

and then he can obtain the value of the initial problem (MP):

$$\frac{\pi^2}{3!} = \frac{\pi^2}{6} = 1 + \frac{1}{4} + \frac{1}{9} + \frac{1}{16} + \ldots + \frac{1}{n^2}$$

It is worth noting here that this isn't in a strict sense the solution of the Mengoli's problem: it is a plausible value suggested by an analogical inference. Euler simply shows how by the soundness of the combined analogical hypotheses (I) and (II) the truth of the assertion $1 + \frac{1}{4} + \frac{1}{9} + \frac{1}{16} + \ldots + \frac{1}{n^2} = \frac{\pi^2}{6}$ would follow.

But the truth of the hypotheses (I) and (II) is not granted at all: it is just another problem to solve that depends, in turn, on other hypotheses to prove. In fact, (I) and (II) are based on inferences which, form a strict logical point of view, are incorrect and bold passages. In particular, in order to prove them, it is necessary to prove firstly the correctness of the passage from (I) finite to infinite, and then (II) from algebra to trigonometry. Such a proof obviously could not, and cannot, be given in general: (I) and (II) are 'dangerous' inferential steps, which can easily lead to false propositions or contradictions (1=0).

Let us consider for simplicity the hypothesis (I), i.e. the analogical transfer of properties from finite sums (i.e. the *source*) into infinite sums (i.e. the *target*). Let us consider the following *infinite* sum:

$$C = 1 - 1 + 1 - 1 + 1 - 1 + \ldots$$

Let apply to C the associative property (according to which the way chosen to arrange the terms of finite sums doesn't affect its result), in the following way:

$$C = (1-1) + (1-1) + (1-1) + \ldots$$

We obtain,

$$C = 0+0+0+...+0 = 0.$$

Now let us apply again the associative property to C as follows:

$$C = 1-(1+1)-(1+1)-(1+1)...$$

We have

$$C = 1-((1+1)-(1+1)-(1+1)...)$$

that is,

$$C = 1-((2-2)-(2-2)-(2-2)...)$$
$$C = 1-((0-0)-(0-0)...)$$
$$C = 1-0 = 1$$

therefore *C = 1 = 0*. We obtained a contradiction.

Obviously in the case of *finite* sums, the same process does not lead to falsity. For example, let *c* be

$$c = 1-1+1-1+1$$

Following the first arrangement of its items we obtain:

$$c = (1-1)+(1-1)+1$$

that is

$$c = 0+0+1 = 1$$

and following the second arrangement of its items:

$$c = 1-(1+1)-(1+1)$$

that is

$$c = 1-((1+1)-(1+1))$$
$$c = 1-(2-2)$$
$$c = 1-0 = 1$$

Therefore, the results are identical. We can conclude that «you cannot make on infinite sums (composition), which have to be previously defined, the same operations allowed for finite sums (compositions)»[55]. Their validity, if provable, is therefore only local.

Thus, in order to increase the degree of reliability of his analogical conclusion, Euler decided to adopt the following principles, as Polya[56] explains, to corroborate its inference:

C1) *a conjecture becomes more credible through the verification of any of its new consequences*

C2) *a conjecture becomes more credible if analogous conjectures become more credible*

According to C1, Polya[57] suggests for example the following consequences which are verifiable to corroborate Euler's analogical inference (obviously, the list can be very long):

- Does the analogical conjecture agree with the 'known fact' $\sin(-x) = -\sin x$?
- Does it agree with $\sin(x+\pi) = -\sin x$?
- Does it agree with $\sin x = 2\sin\left(\frac{x}{2}\right)\cos\left(\frac{x}{2}\right)$?

The principle C2 shows how analogical reasoning can be used at the same time as a means to discover new knowledge and at the same time as a means to justify knowledge. In fact, in order to increase the plausibility of its analogically projected value, Euler applied the *same analogical* hypotheses to an *analogous* problem, namely the Leibniz series, that is $1 - \frac{1}{3} + \frac{1}{5} - \frac{1}{7} + \frac{1}{9} - \ldots \frac{(-1)^n}{2n+1}$, succeeding in finding the value that fits the existing knowledge.

## § 4 The limits of analogy: «analogies are bound to break down even if initially fertile»

No one can reasonably doubt the fruitfulness of analogy: «the ubiquity and fruitfulness of analogy in hypothesis formation is so obvious that they hardly need to be extolled»[58]. Nevertheless, the limits of analogy are equally evident.

First of all, analogy can have only a local validity: an analogical hypothesis (e.g. the finite-infinite analogy for sums in Euler's solution of the problem of Mengoli) can be accepted only under accurate and limited conditions[59]. Secondly and above all, by analogy we extend our knowledge in a very restricted way, i.e. only and simply by suggesting that the unknown behaves somehow like the known and familiar. Surely it is always convenient, also in terms of expected utility (whether or not the analogical hypothesis succeeds), because we shall learn something on the domain under investigation. If the analogical hypothesis passes the tests, we shall learn that the source and target domain are indeed similar (either formally or substantially). If it utterly fails, we shall come to know that we need some radically new ideas to deal with the target domain. Nevertheless, it is precisely what is really and radically new that can not be accounted for by analogy: analogical reasoning can not in principle capture radical new properties of a domain, i.e. what is different from any known propriety or structure (*sui generis*).

---

[55] Sossinsky 1999, 64
[56] Polya 1954, I, 22
[57] ivi, 30-31
[58] Bunge 1967, 265
[59] Hardy1908

Furthermore, analogy exhibits dynamical limits: it can start from fruitfulness and end in nonsense. Quantum mechanics is an example of such dynamical limits, in which an initial analogical success becomes a failure: «in particular analogy between quantum systems and classical particles and waves become a stumbling block preventing a consistent interpretation of the theory»[60]. The result is that the double analogy between classic physics and quantum physics has to be abandoned in order to gain a 'real' understanding of quantum mechanics: «if we want to build or learn new theory then we are likely to use analogy as a bridge between the known and the unknown. But as soon as the new theory is on hand it should be subjected to a critical examination with a view to dismounting its heuristics scaffolding and reconstructing the system in a literal way»[61].

Although Bunge's criticisms of analogy is the consequence of a logical empiricist and realist conception of analogy I disagree with (i.e. analogy is an obstacle because it can't provide a literal and objective description of the quantum world), he points out some important limits (both static and dynamical) of analogy, which not only affect both the demonstrative and the non-demonstrative role of analogy, but should also be taken in account every time analogy is used or analysed.

## §5    Justification and discovery by analogy: the multiple case

One of the most interesting aspects of analogy is its *multiple* nature, which Polya was the first to consider[62]: the possibility to produce «analogies in which *more than one source analog* is used to reason about a target analog»[63]. Multiple analogies have important methodological and philosophical consequences, in particular on the tenability of the accounts of discovery and justification as independent phases of inquiry. In fact, in multiple analogical reasoning, not only analogy shows up in both contexts (discovery and justification), but more than that, the very *same* analogy can be used *at the same time* both as a means of discovery and as a means of justification. Let us consider an example from archaeology, given by Talaly[64].

---

[60] ivi, 265
[61] ivi, 282
[62] Polya 1957
[63] Shelly 2003, 3-4

The problem to solve is (P): to find the *function* of 18 unusual clay fragments found in 5 sites belonging to Neolithic Peloponnese.

In order to solve it, Talaly used multiple analogy both to generate an hypothesis and then to justify it. In short:

    the target of the problem (the 18 fragments) has some properties/pieces of evidence. In fact, they show certain characteristics:

        (a) fragments seem replica of individual, female legs;
        (b) each leg was once half of a pair of legs;
        (c) these pair of legs were manufactured in order to be broken apart.

    Talaly searches for analogs sources agreeing with the target under some of the properties (a)-(b)-(c). Thus, he is able to find three sources:

        i)     Herodotus symbolon;
        ii)    Euripides symbolon;
        iii)   Romans tesserea hospitalis;

The source (i) speaks for (a), while the sources (ii) and (iii) speak for (b):

    (a)    object serves as a contractual device symbolizing an agreement, obligation, friendship or a common bond
    (b)    object serves as an identifying token between individuals or group, symbolizing an agreement, obligation, friendship or a common bond

So Talaly is able to formulate a first analogical hypothesis (H): the target abjects serve or as (a) or as (b).

    At this point, Tataly searches for other analogs sources in order to try to refine the hypothesis and to test and justify it with existing knowledge.

        Thus, he is able to find three sources:
        iv)    Japanese warifu;
        v)     Ancient chine bronze figure;
        vi)    American mafia bill;

Again, the sources (iv) and (vi) speak for (a), while the source (v) for (b). Thus, on one hand they refine and on the other they corroborate and increase the degree of justification of the hypothesis H.

Thus, data used in the discovery and formations of hypothesis are also used to justify the fact that the hypothesis is accepted or rejected (i.e. to increase *the confidence in the hypothesis*). Consequently, we have that «or this practice must be rejected as illegitimate, or construal of discovery and justification as independent phases if inquiry must be rejected»[65].

---

[64] Talalay 1987
[65] Shelley 2003, 134

The untenability of the distinction between the phases of discovery and justification is not new in the literature: it has been explicitly considered by Wylie[66] and independently refined by Cellucci[67]. In particular, Cellucci claims that «discovery is not placed in a preliminary phase of mathematical inquiry, followed by another one, the justification, but covers its entire course»[68], and it is necessary to «acknowledge that the methods of discovery are at the same time also methods of justification»[69].

Cellucci's conception is more in general the result of a new approach to logic and mathematics (i.e. the *heuristic conception*), which is based on open systems and a revision of the analytic method. The pattern of multiple analogies and the cases of the non-demonstrative role of analogy in § 3 well fit Cellucci's heuristic conception and suggest another idea of analogy, which can be defined as the *heuristics conception of analogy*. This conception takes in into account analogy from a completely new point of view, according to which the main aim of a theory of analogy is not the solution of LPA, but the development and the refinement of methods to generate explicative hypotheses to solve the problem under consideration (i.e. an heuristics purpose).

Moreover, the heuristic conception of analogy is based on open systems theory and revised analytic method. According to this conception, analogy operates in an open conceptual system according to the analytic method, in the following bottom-up way:

1) the process, just like in Euler's solution of the problem of Mengoli, stars with a problem to solve and the analytic method allows the recognition of the conditions of solvability of the problem

2) the analogy, by multiple and iterated interactions with the corpus of existing knowledge (analogy is time-sensitive and context-dependent), looks for hypotheses – or a combination of them - to satisfy at least some of the conditions of solvability of the problem

3) the hypotheses so generated are tested within the existing knowledge

4) the hypotheses depend on other hypotheses which have to be proved (or corroborated), refined or specified and which, in turn, depend on other hypotheses and so on (potentially *ad infinitum*).

Thus, looking at analogy under a very different light, the heuristic conception rejects the standard account of analogy and is the basis for the development of an alternative and more sophisticated process of analogical transfer.

---

[66] Wylie 1985
[67] Cellucci 1998 and 2002
[68] Cellucci 2002, 146
[69] *idib.*

As a matter of fact, the inductive and structural accounts of analogy are expression of a *foundationalist* approach, whose purpose is to search for and offer a criterion to make an analogical transfer well-founded, that is at some extent able to solve LPA. As a consequence, they try to reduce each time the analogy to more primitive notions (e.g. similarity or determinations), whose justification, as I argued, is problematic and relies in turn on controversial hypotheses.

By contrast, the heuristics approach, absorbing the paradox of inference, relies on the explicit idea that analogy can't be intrinsically well-founded at all and that solution of LPA is not, and cannot be, the main aim of a theory of analogy (at the limit, it is not a problem at all). As a consequence, analogy is conceived as a dangerous but fruitful inference, whose richness (the capability of generating explicative hypotheses by) relies on its intrinsic non-soundness.

## Bibliographical references